\documentclass[11pt]{article}
\def\a{\alpha}

\def\b{\beta}
\def\bw{\mbox{$\bigwedge$}}

\def\bwsR{\bw^s\R^{m*}}
\def\bwsdashR{\bw^{s'}\R^{m*}}
\def\chibar{\overline{\chi}}
\def\d{\delta}
\def\D{\Delta}
\def\dh{d_{\mathrm{h}}}
\def\dv{d_{\mathrm{v}}}
\def\dT{d_{\mathrm{T}}}

\def\F{\mathcal{F}}

\def\Fmak{\F_{(m)}^{k-1}}

\def\Fmk{\F_{(m)}^k}
\def\Fmka{\F_{(m)}^{k+1}}
\def\Fmks{\F_{(m)}^{k+s}}
\def\Fml{\F_{(m)}^l}
\def\g{\gamma}
\def\GL{\mathrm{GL}}
\def\H{\mathcal{H}}

\def\id{\mathrm{id}}
\def\im{\mathrm{im}}
\def\iT{i_{\mathrm{T}}}

\def\jkg{j^k\g}

\def\jkog{j^k_0\g}
\def\jkaog{j^{k+1}_0\g}
\def\l{\lambda}
\def\L{\Lambda}
\def\N{\mathbf{N}}

\def\Oaa{\Omega^{1,1}}

\def\Oam{\Omega^{1,m}}
\def\Oama{\Omega^{1,m-1}}
\def\Oasa{\Omega^{1,s+1}}
\def\Oao{\Omega^{1,0}}
\def\Oba{\Omega^{2,1}}
\def\Obar{\overline{\Omega}}
\def\Obb{\Omega^{2,2}}
\def\Obm{\Omega^{2,m}}
\def\Obma{\Omega^{2,m-1}}
\def\Obo{\Omega^{2,0}}
\def\Omo{\Omega^{m,0}}
\def\Ooa{\Omega^{0,1}}
\def\Ooabar{\Obar^{0,1}}

\def\Oom{\Omega^{0,m}}
\def\Ooma{\Omega^{0,m-1}}
\def\Oomabar{\Obar^{0,m-1}}
\def\Oombar{\Obar^{0,m}}
\def\Ooo{\Omega^{0,0}}
\def\Ooobar{\Obar^{0,0}}
\def\Oos{\Omega^{0,s}}
\def\Oosa{\Omega^{0,s+1}}
\def\Oosabar{\Obar^{0,s+1}}
\def\OosAbar{\Obar^{0,s-1}}
\def\Oosbar{\Obar^{0,s}}
\def\Or{\Omega^r}
\def\Ordash{\Omega^{r'}}
\def\Ora{\Omega^{r,1}}
\def\OrA{\Omega^{r-1}}
\def\OrAo{\Omega^{r-1,0}}
\def\OrAs{\Omega^{r-1,s}}
\def\Orasa{\Omega^{r+1,s+1}}
\def\OrasA{\Omega^{r+1,s-1}}
\def\OrasB{\Omega^{r+1,s-2}}
\def\OrAsa{\Omega^{r-1,s+1}}
\def\Orbs{\Omega^{r+2,s}}
\def\OrAm{\Omega^{r-1,m}}
\def\Orm{\Omega^{r,m}}
\def\Orma{\Omega^{r,m-1}}
\def\Oro{\Omega^{r,0}}
\def\Ors{\Omega^{r,s}}
\def\Oraa{\Omega^{r+1,1}}
\def\Oras{\Omega^{r+1,s}}
\def\Orsa{\Omega^{r,s+1}}
\def\OrsA{\Omega^{r,s-1}}
\def\p{\partial}

\def\pd#1#2{\frac{\p #1}{\p #2}}
\def\R{\mathbf{R}}

\def\T{\mathbf{T}}
\def\taumbakkE{\tau_{(m)}^{(r+1)(k+1)-1,k}}
\def\taumbakbkE{\tau_{(m)}^{(r+1)(k+1)-1,(r+1)k}}
\def\taumbkkE{\tau_{(m)}^{(r+1)k,k}}
\def\taumkE{\tau_{(m)}^k}
\def\taumklE{\tau_{(m)}^{k,l}}
\def\taumkakE{\tau_{(m)}^{k+1,k}}
\def\taumkkaE{\tau_{(m)}^{k,k-1}}

\def\ve{\varepsilon}
\def\vf#1{\frac{\p}{\p #1}}
\def\vt{\vartheta}

\def\Xibar{\overline{\Xi}}

\def\hook{\kern 3pt \vrule height 0pt depth 0.4pt width 3pt
\vrule height 5pt depth 0.4pt\kern 3pt}

\def\bk{\rule{0em}{1.6ex}}
\def\blob{\bk\hfill\rule{0.5em}{1.5ex}}

\newcommand{\art}[6]{#1: #2 {\it #3\/} {\bf #4} (#5) #6}
\newcommand{\preprint}[4]{#1: #2 {\it Preprint\/} #3 (#4)}

\newcommand{\inbook}[6]{#1: #2 {\it In:\ #3\/} (#4, #5) #6}

\begin{document}

\title{Homogeneous variational complexes and bicomplexes}
\author{D. J. Saunders\\
Department of Algebra and Geometry\\
Palacky Unversity\\
779 00 Olomouc, Czech Republic\footnote{Address for
correspondence: 30 Little Horwood Road, Great Horwood,
Milton Keynes, MK17 0QE, UK}\\
e-mail david@symplectic.demon.co.uk}

\maketitle

\begin{abstract}\noindent
We present a family of complexes playing the same r\^{o}le, for homogeneous
variational problems, that the horizontal parts of the variational bicomplex
play for variational problems on a fibred manifold. We show that, modulo
certain pullbacks, each of these complexes (apart from the first one) is
globally exact. All the complexes may be embedded in bicomplexes, and we show
that, again modulo pullbacks, the latter are locally exact. The edge sequence
is an important part of such a bicomplex, and may be used for the study of
homogeneous variational problems.
\\[1ex]
\end{abstract}

\section{Introduction}

The {\em variational bicomplex}\/, introduced around 30 years ago, plays an
important r\^{o}le in the calculus of variations on a fibration
$\pi : E \to M$. The spaces in this bicomplex are spaces of differential forms
on the infinite jet bundle of the fibration. Each form can be decomposed into a
number of components; this decomposition is generated by the decomposition of a
1-form into horizontal and contact components. The exterior derivative $d$ may
thus be written as the sum of two anticommuting parts, $\dh$ and $\dv$, and
both operators are locally exact. The local exactness of $\dv$ essentially
mirrors that of $d$ on the fibres. However, the local exactness of $\dh$ is
significantly harder to demonstrate, although several proofs are available: see,
for example,~\cite{AndDuch, Tak, Tul, Vin1, Vin2}.

By its nature, the variational bicomplex contains no information about the
order of the forms involved, although it is also possible to define a version of
the bicomplex on finite-order jet bundles (see, for example,~\cite{Vit}).

The {\em finite-order variational complexes} are defined in a somewhat
different way. Here, attention is restricted to a particular $k$-jet manifold.
The first observation is that every contact 1-form on this manifold is
horizontal over the $(k - 1)$-jets. An important result (the $C\Omega$
hypothesis) is that every contact form may be expressed in terms of these
contact 1-forms and their exterior derivatives. Thus a complex may be obtained
from the spaces of forms on the $k$-jet manifold by taking quotients with
respect to the contact forms and their exterior derivatives. The exterior
derivative $d$ passes to this quotient, and is locally exact~\cite{KrupSeq}.

The finite-order variational complex can also be defined for the manifold of
$k$-th order contact elements of a manifold $E$. However the variational
bicomplex makes essential use of the fibration over the base manifold $M$.

In this paper, we consider the possiblity of defining similar complexes
for homogeneous variational problems. These are problems defined on $m$-frame
bundles, that is, bundles of regular $m$-velocities (see, for
example,~\cite{CS1, CS2}). The prototype for such problems is Finsler geometry,
where the problem is defined on the slit tangent bundle of $E$ (in other words,
the first-order 1-frame bundle). Any variational problem on a jet bundle of a
fibration induces a homogeneous variational problem on (an open submanifold of)
a frame bundle and studying the problem in this context can sometimes give
important insights.

Our definition is given for finite-order $m$-frame bundles. The terms
in each complex contain vector-valued forms rather than scalar forms, and we
show that, modulo certain pull-backs, each complex is exact. We also embed
these complexes in suitable bicomplexes, where the vertical differential $\dv$
is replaced by the ordinary exterior derivative $d$. The edge sequence of this
bicomplex has some similarities with finite-order variational complex on a
jet manifold.

We start, therefore, with some preliminary remarks about the frame bundles and
their properties, and then introduce the spaces of each complex and the
coboundary operator $\dT$. The main result of the paper is the proof that, for
all except the first complex, this operator is globally exact (modulo
pullbacks) and, indeed, that there is a canonically-defined `homotopy operator'
--- the quotation marks indicate that it takes its values in a pull-back of
the domain of $\dT$, so that the term {\em pseudo-homotopy operator}\/ might be
more appropriate. The following section introduces the bicomplexes and
demonstrates (local) exactness of the first complex, again modulo pullbacks,
while the final section shows how certain terms of the edge sequence containing
equivalence classes of vector-valued forms may be mapped globally to spaces of
representative forms, and considers the relationship of this sequence to
problems in the calculus of variations. A subsequent paper~\cite{SauFund} will
give a preliminary report on a project to apply this theory to find, for a
homogeneous Lagrangian, a scalar form which is closed precisely when the
Lagrangian is null:\ this corresponds to the `fundamental Lepage equivalent'
known only for first-order Lagrangians in the case of fibred manifolds.

\section{Preliminary remarks}

The manifolds studied in the context of homogeneous variational problems are
the (higher-order) $m$-frame bundles (see, for example,~\cite{CS2}). Given a
smooth manifold $E$ with $\dim E = n$, consider the $k$-th order $m$-velocities
in $E$. These are the $k$-jets at $0 \in \R^m$ of maps $\g$ from a
neighbourhood of $0$, to $E$. If a map $\g$ is an immersion then its $k$-jet
$\jkog$ is called a $k$-th order $m$-frame. The set of all such $\jkog$ is
denoted $\Fmk E$. It is a smooth manifold, and a bundle over $E$ with
projection $\taumkE : \Fmk E \to E$. There are also projections
$\taumklE : \Fmk E \to \Fml E$ where $l < k$. If $(u^\a)$ are local coordinates
on $U \subset E$, where $1 \leq \a \leq n$, then
\[
(u^\a, u^\a_i, u^\a_{(ij)}, \ldots, u^\a_{(i_1 i_2 \cdots i_k)})
\]
are local coordinates on $(\taumkE)^{-1}(U) \subset \Fmk E$, where
$1 \leq i, j, \ldots \leq m$ and the parentheses $(ij)$ indicate symmetrization.
In multi-index notation, these coordinates may be written as $(u^\a_I)$ where
$I \in \N^m$ and where $0 \leq |I| \leq k$. Note that the two notations often
give rise to different constant multiples.

A basic example of this construction arises when $k = 1$. In this case the
bundle of first-order $m$-velocities may be identified with the Whitney sum
bundle $\bigoplus^m TE \to E$, and the bundle of first-order $m$-frames may be
identified with the set of $m$-tuples $(\xi_1, \xi_2, \ldots, \xi_m)$ where the
vectors $\xi_i$ are linearly independent. For this example, there is another
interpretation of the bundle of $m$-velocities as a tensor product bundle
$TE \otimes \R^{m*}$. If $(\xi_1, \xi_2, \ldots, \xi_m)$ is a first-order
$m$-velocity (where $\xi_i \in T_p E)$ then the corresponding tensor is
$\xi_i \otimes t^i$, where $(e_i)$ is the standard basis of $\R^m$ and
$(t^i)$ is the dual basis of $\R^{m*}$.

Frame bundles are closely related to bundles of contact elements, and to jet
bundles. The bundle $J^k (E, m)$ of $m$-dimensional $k$-contact elements on $E$,
also known as the bundle of $k$-jets of immersed $m$-dimensional submanifolds, or
the bundle of $k$-th order Grassmannians, is a quotient of the frame bundle
$\Fmk E$. The group of diffeomorphisms of $\R^m$ preserving the origin acts on
elements of $\Fmk E$ by composition with the immersion defining the frame, and
this action factors through to an action of the $k$-th order jet group; if
$k = 1$ then this group is $\GL(m)$. The quotient of $\Fmk E$ by this action
of the jet group is $J^k (E, m)$. There is also an oriented version of this
quotient construction, given by taking the subgroup of orientation-preserving
diffeomorphisms.

Now suppose that there is a fibration $\pi : E \to M$ over some $m$-dimensional
manifold $M$. The image of any local section $\phi$ of $\pi$ is an
$m$-dimensional submanifold of $E$, and so any jet $j^k_p\phi \in J^k \pi$ is
an element of $J^k (E, m)$; in this way $J^k\pi$ is an open submanifold of
$J^k (E, m)$. It is not the whole of $J^k (E, m)$, because submanifolds that
are not transverse to the fibration do not have jets in $J^k \pi$.

In the special case where $M = \R^m$, each local section of $\pi$ gives rise by
translation to a map from a neighbourhood of the origin to $E$, and so the jet
$j^k_p\phi$ gives rise to an $m$-velocity; it is in fact an $m$-frame because a
local section has rank $m$. But only those maps $\R^m \to E$ that are (local)
sections of $\pi$ have velocities corresponding to jets in $J^k \pi$, and so
in this particular special case $J^k \pi$ becomes a closed submanifold of
$\Fmk E$.

A choice of adapted coordinates $(x^i, y^a)$ on $E$ gives a local
identification of $M$ with $\R^m$, and so locally $J^k \pi$ may be considered
as a closed submanifold of $\Fmk E$ (but not in an invariant way: the
identification depends on the chart). Put $u^i = x^i$ for $i = 1, \ldots, m$
and $u^{m + a} = y^a$ for $a = 1, \ldots, n - m$; then $J^k \pi$ is the
submanifold of $\Fmk E$ gven by $u^i_j = \d^i_j$, $u^i_I = 0$ for $|I| > 1$.

There are a number of objects associated with frame bundles, and corresponding
formul\ae\ relating them. Most of the objects are generalizations of those to be
found on tangent bundles. We shall need to consider, in particular, the total
derivatives and the vertical endomorphisms.

First, the total derivatives. These arise because a prolonged map
$\jkg : \R^m \to \Fmk E$ defines both a point $\jkaog$ in the $(k+1)$-th frame
bundle $\Fmka E$, and a first-order $m$-velocity $j^1_0(\jkg)$ at
$\jkog \in\Fmk E$. This relationship gives an embedding
$\T : \Fmka E \to T(\Fmk E) \otimes \R^{m*}$ called the {\em $(k + 1)$-th
order total derivative for $m$-frames}. The $i$-th component of $\T$ is the
contraction $\T_i = \langle \T, e_i \rangle$, and is a vector field along
the projection $\taumkakE : \Fmka E \to \Fmk E$; in coordinates it is
\[
\T_i = \sum_{|I| = 0}^k u^\a_{I+1_i} \vf{u^\a_I} \, .
\]
We shall write $d_i$ for the action of $\T_i$ on the functions on $\Fmk E$,
giving functions on $\Fmka E$, and of course this may be extended to a map of
$r$-forms $d_i : \Or \Fmk E \to \Or \Fmka E$. We also write
$i_i : \Or \Fmk E \to \OrA \Fmka E$ for the corresponding contraction.

The `vertical endomorphisms' may be defined in the following way. Suppose that
$\chi : \R^m \times \R^m \to E$ is a map, and put $\chi_y(x) = \chi(x, y)$;
suppose that each $k$-th order $m$-velocity $j^k \chi_y$ is a $k$-th order
$m$-frame, and so an element of $\Fmk E$. The 1-jet at zero of the map
$\R^m \to \Fmk E$, $y \mapsto j^k \chi_y$, is a 1st-order velocity on
$\Fmk E$ at $j^k \chi_0$, and so is an element of $T\Fmk E \otimes \R^{m*}$;
any such velocity may be represented by a map $\chi$ in this way.

Given $\chi$, define a new map $\chibar : \R^m \times \R \to E$ by
$\chibar(x, t) = \chi(x, tx)$, and write $\chibar_t(x) = \chibar(x, t)$. The
map $\R \to \Fmka E$, $t \mapsto j^{k + 1} \chibar_t$ is a curve in $\Fmka E$
defining a tangent vector in $\Fmka E$.

Given a 1st-order velocity in $T\Fmk E \otimes \R^{m*}$ at $a \in \Fmk E$ and a
point $b \in \Fmka E$ projecting to $a$, we may represent the 1-velocity by a
map $\chi$ as described above, choosing $\chi$ so that $j^{k + 1} \chi_0 = b$.
We may then construct the corresponding map $\chibar$, noting that the curve
$t \mapsto j^{k + 1} \chibar_t$ passes through $b$ at zero. We thus obtain a
tangent vector at $b$, and it may be shown that this is independent of the
choice of the representative map $\chi$. This tangent vector is the vertical
lift of the velocity to the point $b$.

We may now map a tangent vector $\xi$ at $b \in \Fmka E$ to another tangent
vector at the same point. First we project $\xi$ to $a \in \Fmk E$; then we
construct the velocity $(0, \ldots, \taumkakE(\xi), \ldots, 0)$ with the
projected tangent vector in the $i$-th place; and finally we take the vertical
lift to $b$. We denote the resulting tangent vector by $S^i(\xi)$. The
operations are all linear, and so $S^i$ is a tensor field of type $(1, 1)$ on
$\Fmka E$. These are the vertical endomorphisms of the frame bundles. In
coordinates, they are
\[
S^i = \sum_{|I|=0}^k (I(i) + 1) \vf{u^\a_{I + 1_i}} \otimes du^\a_I \, .
\]

\section{Homogeneous variational complexes}

The variational bicomplex is constricted from spaces of (scalar) differental
forms, and the finite-order variational sequences are constructed from spaces
of equivalence class of scalar forms. By contrast, the homogeneous variational
complexes that we shall describe below are constructed from spaces of {\em
vector-valued forms}\/: to be precise, we consider the vector spaces $\bwsR$ of
multilinear forms on $\R^m$, and study the differential forms on $\Fmk E$
taking their values in these vector spaces. We shall denote the set of such
$r$-forms by $\Ors_k = \Or_k \otimes \bwsR$, where $\Or_k = \Or\Fmk E$ is the
space of scalar $r$-forms on $\Fmk E$. In coordinates, an
element of $\Ors_k$ would be written
\[
\Phi = \Phi^{I_1 \cdots I_r}_{a_1 \cdots a_r i_1 \cdots i_s}
\left( du^{a_1}_{I_1} \wedge \ldots \wedge du^{a_r}_{I_r} \right)
\otimes \left( dt^{i_1} \wedge \ldots \wedge dt^{i_s} \right) \, ,
\]
where $t^i$ are the standard linear coordinate functions on $\R^m$ and we write
$dt^i$ rather than $t^i$, using the identification between the constant 1-form
$dt^i$ and the element $t^i\in \R^{m*}$ of the dual space.

The wedge products on the two components of these vector-valued forms induce a
wedge product on the direct sum $\bigoplus_{r,s} \Ors_k$: for decomposable
elements,
\[
(\theta \otimes w) \wedge (\theta' \otimes w')
= (\theta \wedge \theta') \otimes (w \wedge w')
\]
where $\theta \in \Or_k$, $\theta' \in \Ordash_k$, $w \in \bwsR$,
$w' \in \bwsdashR$. The fact that, for $\l \in \R$, $\l \not= 0$,
\[
(\l \theta) \otimes (\l^{-1} w) = \theta \otimes w
\]
causes no ambiguity in this definition.

The exterior derivative on the first component induces a derivation
$d : \Ors_k \to \Oras_k$: for decomposable elements
\[
d(\theta \otimes w) = d\theta \otimes w \, ,
\]
and this is extended by linearity. Again there is no ambiguity in this
definition.

Finally, the vertical endomorphisms $S^i$ may be combined in this context to
produce a map $S : \Ors_k \to \OrsA_k$ by writing
\[
S(\theta \otimes w) = S^i\theta \otimes w_i \, ,
\]
where $w_i$ denotes the contraction of $\p / \p t^i$ with $w$.

We now introduce two new operations on these vector-valued forms, both derived
from the total derivative operators on scalar forms. These are the maps
$\iT : \Ors_k \to \OrAsa_{k + 1}$ and $\dT : \Ors_k \to \Orsa_{k + 1}$,
defined for decomposable elements by
\begin{eqnarray*}
\iT(\theta \otimes w) & = & i_i \theta \otimes (dt^i \wedge w) \\
\dT(\theta \otimes w) & = & d_i \theta \otimes (dt^i \wedge w)
\end{eqnarray*}
and extended by linearity. As the $i_i$ and $d_i$ are derivations, it follows
that $\iT$ and $\dT$ are derivations. Note that, in general, both $\iT$ and
$\dT$ increase the order of a vector-valued form by one. We call $\dT$ the
{\em total exterior derivative}.

We can now build the homogenous variational complexes from the spaces of
vector-valued forms and the maps $\dT$. It is immediate that $\dT^2 = 0$,
because
\[
\dT^2(\theta \otimes w) = d_i d_j \theta \otimes (dt^i \wedge dt^j\wedge w)
= 0
\]
as $d_i d_j$ is symmetric in $i, j$ whereas $dt^i \wedge dt^j$ is
skew-symmetric. Thus for each $k$ and each $r$ satisfying
$1 \leq r \leq \dim \Fmk E$ we have a complex
\[
0 \to \Oro_k \to \Ora_{k+1} \to \ldots \to \Ors_{k+s} \to \Orsa_{k+s+1}
\to \hspace{4em}
\]
\\[-6ex]
\[
\hspace{4em} \to \ldots \to \Orm_{k+m} \to \Orm_{k+m} / \dT(\Orma_{k+m-1})
\to 0 \, ,
\]
and for each $k$ when $r = 0$ a complex
\[
0 \to \R \to \Ooo_k \to \Ooa_{k+1} \to \ldots \to \Oos_{k+s} \to \Oosa_{k+s+1}
\to \hspace{4em}
\]
\\[-6ex]
\[
\hspace{4em} \to \ldots \to \Oom_{k+m} \to \Oom_{k+m} / \dT(\Ooma_{k+m-1})
\to 0 \, .
\]
We shall, in fact, replace the latter complex by
\[
0 \to \Ooobar_k \to \Ooabar_{k+1} \to \ldots \to \Oosbar_{k+s} \to
\Oosabar_{k+s+1} \to \hspace{4em}
\]
\\[-6ex]
\[
\hspace{4em} \to \ldots \to \Oombar_{k+m} \to
\Oombar_{k+m} / \dT(\Oomabar_{k+m-1}) \to 0
\]
where $\Oosbar_{k+s} = \Oos_{k+s} / \bw^s \R^{m*}$, taking quotients by the
constant vector-valued functions. The reason for this will be given when we
consider exactness.

There are also some shorter complexes arising because we are considering
finite-order $m$-frame bundles. We have
\[
0 \to \Ors_0 \to \Orsa_1 \to \ldots \to \Orm_{m-s}
\to \Orm_{m-s} / \dT(\Orma_{m-s-1}) \to 0
\]
for $r \geq 1$ and
\[
0 \to \Oosbar_0 \to \Oosabar_1 \to \ldots \to \Oombar_{m-s}
\to \Oombar_{m-s} / \dT(\Oomabar_{m-s-1}) \to 0
\]
arising when we consider zeroth-order forms taking their values in
$\bw^s \R^{m*}$ for $s \geq 1$, and
\[
\Ors_k \to \Orsa_{k+1} \to \ldots \to \Orm_{k+m-s}
\to \Orm_{k+m-s} / \dT(\Orma_{k+m-s-1}) \to 0
\]
arising when $r > \dim \Fmak E$ (so that $\OrsA_{k-1} = 0$). We shall call each
of the complexes described above a {\em homogeneous variational complex}.

\section{Exactness}

The main result of this paper concerns the exactness of the homogeneous
variational complexes described above. A first remark is that it is, in general,
not possible for all of these to be exact in the usual sense, even locally. To
see why, let $E = \R^n$ and let $u^\a$ be the standard (global) coordinate
functions on $E$. Let $\Phi$ be the non-zero element of $\Obb_1$ given in these
coordinates by
\[
\Phi = \d_{\a\b} du^\a_i \wedge du^\b_j \otimes dt^i \wedge dt^j,
\]
so that
\[
\dT \Phi = \d_{\a\b} \left(
du^\a_{ik} \wedge du^\b_j + du^\a_i \wedge du^\b_{jk} \right)
\otimes dt^k \wedge dt^i \wedge dt^j = 0.
\]
But for any $\Psi \in \Oba_0$ we must have
\[
\Psi = \Psi_{\a\b j} du^\a \wedge du^\b \otimes dt^j
\]
with $\Psi_{\a\b j} + \Psi_{\b\a j} = 0$, so that
\[
\dT \Psi = \left( (d_i \Psi_{\a\b j}) du^\a \wedge du^\b
+ 2\Psi_{\a\b j} du^\a_i \wedge du^\b \right) \otimes dt^i \wedge dt^j
\not= \Phi.
\]

It is clear that this apparent setback arises because total derivatives
increase the order of forms: indeed, exactly the same phenomenon arises in the
affine case, and this is why the variational bicomplex is usually defined on
the infinite jet bundle. We could adopt a similar procedure here, but again we
would lose all information about the order of the form. Instead, we shall look
at particular {\em finite} pullbacks, and we shall see that, modulo these
pull-backs, the homogeneous variational complexes are globally exact. Indeed,
the result is even stronger than this: we are able to construct a canonical
`pseudo-homotopy operator', so that if $\dT \Phi = 0$ then we can make a
canonical choice of $\Psi$ such that $\dT \Psi = (\taumbkkE)^* \Phi$.

\textbf{Definition}
{\em
If $r \geq 1$, the operator $P : \Ors_k \to \Orsa_{(r+1)k-1}$ is defined by
\begin{eqnarray*}
P(\Phi) & = & P^j_{(s)}(\phi_{i_1 \cdots i_s}) \otimes
\left\{ \vf{t^j} \hook \left( dt^{i_1} \wedge \ldots
\wedge dt^{i_s} \right) \right\} \\
& = & s \, P^j_{(s)}(\phi_{j i_2 \cdots i_s})
\otimes \left( dt^{i_2} \wedge \ldots \wedge dt^{i_s} \right)
\end{eqnarray*}
where $\Phi = \phi_{i_1 \cdots i_s} \otimes \left( dt^{i_1} \wedge \ldots
\wedge dt^{i_s} \right)$, the $\phi_{i_1 \cdots i_s}$ are scalar $r$-forms,
completely skew-symmetric in the indices $i_1, \ldots, i_s$, and $P^j_{(s)}$ is
the differential operator on scalar $r$-forms defined by
\[
P^j_{(s)} = \sum_{|J| = 0}^{rk-1}
\frac{(-1)^{|J|}(m - s)! |J|!}{r^{|J|+1} (m - s + |J| + 1)! J!}
d_J S^{J + 1_j} \, .
\]
}

Note that the scalar forms $\phi_{i_1 \cdots i_s}$ are globally well-defined.
Indeed, $\Phi$ takes its values in the vector space $\bwsR$, and
$\phi_{i_1 \cdots i_s}$ is obtained by contracting the image of $\Phi$ with
the $s$-vector $e_{i_1} \wedge \ldots \wedge e_{i_s}$ where $\{e_i\}$ is the
standard basis of $\R^m$.

\textbf{Theorem}
{\em
The operator $P$ is a pseudo-homotopy operator for $\dT$ for the case
$r \geq 1$, in that
\[
\dT P + P \dT = \id
\]
modulo pullbacks; thus the homogeneous variational complexes are globally
exact, again modulo pullbacks. This formula is also valid for the case
$s = 0$ where, by default, $P^j_{(0)} = 0$, and for the case $s>0, k=0$ where
$P^j_{(s)} = 0$ explicitly.
}

A remark is needed about the spaces on which the homotopy formula is defined.
In principle, if $\Phi \in \Ors_k$ then $\dT P \Phi \in \Ors_{(r+1)k}$ and
$P \dT \Phi \in \Ors_{(r+1)(k+1)-1}$, so that the formula would be an equation
in $\Ors_{(r+1)(k+1)-1}$; to be precise it would be
\[
\left( \taumbakbkE \right)^* \dT P + P \dT
= \left( \taumbakkE \right)^* \, .
\]
But in fact it turns out that $P \dT \Phi$ is always projectable to
$\Ors_{(r+1)k}$, and indeed this is essential for the success of the proof; in
the particular case of interest, where $\dT \Phi = 0$, both sides of the
formula are  projectable further to $\Ors_k$. Note however that, for the case
$r > \dim \Fmak E$ where the complex starts
$\Ors_k \to \Orsa_{k+1} \to \ldots$, the first map is {\em not}\/ injective,
and indeed the codomain of $P : \Ors_k \to \OrsA_{2k-1}$ is not trivial.

We also mention that, for certain classes of forms, $P\Psi$ may be projectable
to a lower-order frame bundle; this turns out to be important for variational
problems (see~\cite{SauFund}).

In the proofs below we shall normally omit the pullback maps from the formul\ae\
where there is no chance of confusion.

The proof of the theorem is based upon the relationship between the total
derivative operators $d_j$ and the vertical endomorphisms $S^i$. The
fundamental lemma is this one.

\textbf{Lemma 1}
{\em
If $\theta$ is a scalar 1-form on $\Fmk E$ then
\[
S^i d_j \theta - d_j S^i \theta = \d^i_j \, (\taumkakE)^* \theta \, .
\]
}

\textbf{Proof}
This result is mentioned in~\cite{CS1}; we include a proof for completeness.
We use local coordinates. Let
\[
\theta = \sum_{|J|=0}^k \theta^J_\a du^\a_J \, ,
\]
and note that we may rewrite $S^i$ as
\[
S^i = \sum_{|I|=0}^{k-1} (I(i) + 1) du^\a_I \otimes \vf{u^\a_{I + 1_i}}
= \sum_{|J|=1}^k J(i) du^\a_{J-1_i} \otimes \vf{u^\a_J}
\]
(the multi-index $J - 1_i$ does not make sense when $J(i) = 0$, but then the
coefficient of the term vanishes so that the formula is still valid.) Thus
\[
d_j \theta = \sum_{|J|=0}^k \left( (d_j \theta^J_\a) du^\a_J
+ \theta^J_\a du^\a_{J+1_j} \right)
\]
so that
\[
S^i d_j \theta = \sum_{|J|=0}^k \left( J(i) (d_j \theta^J_\a) du^\a_{J-1_i}
+ (J + 1_j)(i) \theta^J_\a du^\a_{J+1_j-1_i} \right) \, ;
\]
also
\[
S^i \theta = \sum_{|J|=0}^k J(i) \theta^J_\a du^\a_{J-1_i}
\]
so that
\[
d_j S^i \theta = \sum_{|J|=0}^k J(i) \left(
(d_j \theta^J_\a) du^\a_{J-1_i} + \theta^J_\a du^\a_{J-1_i+1_j} \right) \, .
\]
Subtracting,
\begin{eqnarray*}
S^i d_j \theta - d_j S^i \theta & = & \sum_{|J|=0}^k [
(J + 1_j)(i) - J(i) ] \theta^J_\a du^\a_{J+1_j-1_i} \\
& = & \sum_{|J|=0}^k \d^i_j \theta^J_\a du^\a_{J+1_j-1_i} \\
& = & \d^i_j \theta \, .
\end{eqnarray*}
\blob

\textbf{Lemma 2}
{\em
If $\phi$ is a scalar $r$-form then
\[
S^i d_j \phi - d_j S^i \phi = r \d^i_j \phi \, .
\]
}

\textbf{Proof}
By induction on $r$, using the fact that the commutator of two derivations is
again a derivation.
\blob

We shall now assume that all scalar forms are $r$-forms, and write
\[
S^i d_j - d_j S^i = r \d^i_j
\]
without further comment.

\textbf{Lemma 3}
{\em
For any natural number $p$,
\[
(d_{i_1} \ldots d_{i_p} S^{i_1} \ldots S^{i_p}) d_j
= d_j (d_{i_1} \ldots d_{i_p} S^{i_1} \ldots S^{i_p})
+ rp \, d_j (d_{i_1} \ldots d_{i_{p-1}} S^{i_1} \ldots S^{i_{p-1}}) \, .
\]

}
\textbf{Proof}
By induction on $p$:
\begin{eqnarray*}
(d_{i_1} \ldots d_{i_p} S^{i_1} \ldots S^{i_p}) d_j
& = & (d_{i_1} \ldots d_{i_p} S^{i_1} \ldots S^{i_{p-1}})
(d_j S^{i_p} + r\d^{i_p}_j) \\
& = & [d_j (d_{i_1} \ldots d_{i_p} S^{i_1} \ldots S^{i_p}) \\
& & \phantom{dt} + r(p - 1) d_j d_{i_p} (d_{i_1} \ldots d_{i_{p-2}} S^{i_1}
\ldots S^{i_{p-2}}) S^{i_p}] \\
& & \phantom{dt} + r \, d_j (d_{i_1} \ldots d_{i_{p-1}} S^{i_1} \ldots
S^{i_{p-1}}) \\
& = & d_j (d_{i_1} \ldots d_{i_p} S^{i_1} \ldots S^{i_p}) \\
& & \phantom{dt} + rp \, d_j (d_{i_1} \ldots d_{i_{p-1}} S^{i_1} \ldots
S^{i_{p-1}}) \, ,
\end{eqnarray*}
using Lemma 2, the induction hypothesis, and the commutativity of total
derivatives.
\blob

We shall need to rewrite Lemma~3 using multi-index notation.

\textbf{Lemma 4}
{\em
For any natural number $p$,
\[
\sum_{|I|=p} \frac{|I|!}{I!} d_I S^I d_j
= d_j \left( \sum_{|I|=p} \frac{|I|!}{I!} d_I S^I
+ rp \sum_{|J|=p-1} \frac{|J|!}{J!} d_J S^J \right) \, .
\]
}

\textbf{Proof}
For any multi-index $I$ with
\[
I = 1_{i_1} + \ldots + 1_{i_p} \, ,
\]
the number of {\em distinct} rearrangements of the indices $i_1, \ldots, i_p$
(giving the same multi-index $I$) is $|I|! / I!$ ; this is called the {\em
weight} of $I$.
\blob

We now give a proof of the main theorem, using Lemma~4. As before, we take
$\Phi \in \Ors_k$ with $\Phi = \phi_{i_1 \cdots i_s} \otimes \left( dt^1 \wedge
\ldots \wedge dt^s \right)$; then
\begin{eqnarray*}
\dT P(\Phi) & = & \dT \left\{ s \, P^j_{(s)}(\phi_{j i_2 \cdots i_s})
\otimes \left( dt^{i_2} \wedge \ldots \wedge dt^{i_s} \right) \right\} \\
& = & s \left( d_{i_1} P^j_{(s)}(\phi_{j i_2 \cdots i_s}) \right)
\otimes \left( dt^{i_1} \wedge \ldots \wedge dt^{i_s} \right) \, ,
\end{eqnarray*}
whereas
\begin{eqnarray*}
P( \dT \Phi) & = & P \left\{ (d_j \phi_{i_1 i_2 \cdots i_s})
\otimes \left( dt^j \wedge dt^{i_1} \wedge \ldots \wedge dt^{i_s} \right)
\right\} \\
& = & P^q_{(s+1)} (d_j \phi_{i_1 i_2 \cdots i_s}) \otimes
\left\{ \vf{t^q} \hook \left( dt^j \wedge dt^{i_1} \wedge \ldots
\wedge dt^{i_s} \right) \right\} \\
& = & P^j_{(s+1)} (d_j \phi_{i_1 i_2 \cdots i_s}) \otimes
\left( dt^{i_1} \wedge \ldots \wedge dt^{i_s} \right) \\
& & \phantom{dt} - s P^{i_1}_{(s+1)} (d_j \phi_{i_1 i_2 \cdots i_s}) \otimes
\left( dt^j \wedge dt^{i_2} \wedge \ldots \wedge dt^{i_s} \right) \\
& = & \left( P^j_{(s+1)} (d_j \phi_{i_1 i_2 \cdots i_s})
- s P^j_{(s+1)} (d_{i_1} \phi_{j i_2 \cdots i_s}) \right) \otimes
\left( dt^{i_1} \wedge \ldots \wedge dt^{i_s} \right) \, ,
\end{eqnarray*}
using the definitions of $\dT$ and $P$. So the task is to show that, as scalar
$r$-forms,
\[
s \, d_{i_1} P^j_{(s)}(\phi_{j i_2 \cdots i_s})
+ P^j_{(s+1)} (d_j \phi_{i_1 i_2 \cdots i_s})
- s \, P^j_{(s+1)} (d_{i_1} \phi_{j i_2 \cdots i_s})
\stackrel{?}{=} \phi_{i_1 i_2 \cdots i_s} \, ,
\]
and we carry out this task by examining the operator
\[
s \, d_{i_1} P^j_{(s)} + \d^j_{i_1} P^q_{(s+1)} d_q - s \, P^j_{(s+1)} d_{i_1}
\]
acting on $\phi_{j i_2 \cdots i_s}$. We expand each of the three terms, using
Lemma~2 and Lemma~4 to ensure that all the total derivatives $d_j$ are moved to
the left of all the vertical endomorphisms $S^i$.

First,
\[
s \, d_{i_1} P^j_{(s)} = s \sum_{|J| = 0}^{rk-1}
\frac{(-1)^{|J|}(m - s)! |J|!}{r^{|J|+1} (m - s + |J| + 1)! J!}
d_{J + 1_{i_1}} S^{J + 1_j} \, ;
\]
next,
\begin{eqnarray*}
\d^j_{i_1} P^q_{(s+1)} d_q
& = & \d^j_{i_1} \sum_{|J| = 0}^{r(k+1)-1} \frac{(-1)^{|J|}(m - s - 1)! |J|!}{
r^{|J|+1} (m - s + |J|)! J!} d_J S^{J + 1_q} d_q \\
& = & \d^j_{i_1} \sum_{|J| = 0}^{r(k+1)-1} \frac{(-1)^{|J|}(m - s - 1)! |J|!}{
r^{|J|+1} (m - s + |J|)! J!} d_J S^J (d_q S^q + mr) \\
& = & \d^j_{i_1} \sum_{p = 0}^{r(k+1)-1} \frac{(-1)^p (m - s - 1)!}{r^{p+1}
(m - s + p)!} \left\{ \sum_{|J|=p} \frac{|J|!}{J!} d_{J + 1_q} S^{J + 1_q}
\right. \\
& & \left. \phantom{bigskip} + rp \sum_{|J|=p-1} \frac{|J|!}{J!} d_{J+ 1_q}
S^{J + 1_q} + mr \sum_{|J|=p} \frac{|J|!}{J!} d_J S^J \right\} \\
& = & \d^j_{i_1} \sum_{|J|=0}^{r(k+1)-1} \frac{(-1)^{|J|} (m - s - 1)! |J|!}{
r^{|J|+1} (m - s + |J|)! J!} d_{J + 1_q} S^{J + 1_q} \\
& & \phantom{bigskip} + \d^j_{i_1} \sum_{|J|=0}^{r(k+1)-2}
\frac{(-1)^{|J|+1} (m - s - 1)! (|J| + 1)!}{r^{|J|+1} (m - s + |J| + 1)! J!}
d_{J + 1_q} S^{J + 1_q} \\
& & \phantom{bigskip} + m \, \d^j_{i_1} \sum_{|J|=0}^{r(k+1)-1} \frac{(-1)^{|J|}
(m - s - 1)! |J|!}{r^{|J|} (m - s + |J|)! J!} d_J S^J \, ;
\end{eqnarray*}
and finally, using similar manipulations,
\begin{eqnarray*}
- s \, P^j_{(s+1)} d_{i_1}
& = & - s \sum_{|J| = 0}^{r(k+1)-1} \frac{(-1)^{|J|}(m - s - 1)! |J|!}{
r^{|J|+1}
(m - s + |J|)! J!} d_J S^{J + 1_j} d_{i_1} \\
& = & - s \sum_{|J|=0}^{r(k+1)-1} \frac{(-1)^{|J|} (m - s - 1)! |J|!}{
r^{|J|+1}
(m - s + |J|)! J!} d_{J + 1_{i_1}} S^{J + 1_j} \\
& & \phantom{bigskip} - s \sum_{|J|=0}^{r(k+1)-2}
\frac{(-1)^{|J|+1} (m - s - 1)! (|J| + 1)!}{r^{|J|+1} (m - s + |J| + 1)! J!}
d_{J + 1_{i_1}} S^{J + 1_j} \\
& & \phantom{bigskip} - s \, \d^j_{i_1} \sum_{|J|=0}^{r(k+1)-1} \frac{(-1)^{|J|}
(m - s - 1)! |J|!}{r^{|J|} (m - s + |J|)! J!} d_J S^J \, .
\end{eqnarray*}
So altogether there are seven sums to consider:
\begin{eqnarray}
\makebox[5ex][l]{$s \, d_{i_1} P^j_{(s)} + \d^j_{i_1} P^q_{(s+1)} d_q
- s \, P^j_{(s+1)} d_{i_1}$} \nonumber \\[1ex]
& = & s \sum_{|J| = 0}^{rk-1}
\frac{(-1)^{|J|}(m - s)! |J|!}{r^{|J|+1} (m - s + |J| + 1)! J!}
d_{J + 1_{i_1}} S^{J + 1_j} \\
& & \phantom{skip} + \d^j_{i_1} \sum_{|J|=0}^{rk-1} \frac{(-1)^{|J|}
(m - s - 1)! |J|!}{r^{|J|+1} (m - s + |J|)! J!} d_{J + 1_q} S^{J + 1_q} \\
& & \phantom{skip} + \d^j_{i_1} \sum_{|J|=0}^{rk-1}
\frac{(-1)^{|J|+1} (m - s - 1)! (|J| + 1)!}{r^{|J|+1} (m - s + |J| + 1)! J!}
d_{J + 1_q} S^{J + 1_q} \\
& & \phantom{skip} + m \, \d^j_{i_1} \sum_{|J|=0}^{rk} \frac{(-1)^{|J|}
(m - s - 1)! |J|!}{r^{|J|} (m - s + |J|)! J!} d_J S^J \\
& & \phantom{skip} - s \sum_{|J|=0}^{rk-1} \frac{(-1)^{|J|}
(m - s - 1)! |J|!}{r^{|J|+1} (m - s + |J|)! J!} d_{J + 1_{i_1}} S^{J + 1_j} \\
& & \phantom{skip} - s \sum_{|J|=0}^{rk-1}
\frac{(-1)^{|J|+1} (m - s - 1)! (|J| + 1)!}{r^{|J|+1} (m - s + |J| + 1)! J!}
d_{J + 1_{i_1}} S^{J + 1_j} \\
& & \phantom{skip} - s \, \d^j_{i_1} \sum_{|J|=0}^{rk} \frac{(-1)^{|J|}
(m - s - 1)! |J|!}{r^{|J|} (m - s + |J|)! J!} d_J S^J \, ,
\end{eqnarray}
where we have reduced the maximum summation in sums~(2), (3), (5) and~(6) to
$|J|=rk-1$, and in sums~(4) and~(7) to $|J|=rk$: this is because the scalar
forms $\phi_{j i_2 \cdots i_s}$ in the domain of this operator are defined on
$\Fmk E$, so that $S^J \phi_{j i_2 \cdots i_s} = 0$ whenever $|J| > rk$.

We now show that the seven sums collapse to give $\d^{i_1}_j$.

First, we see that sums~(1), (5) and~(6) cancel when taken together. The range
of summation is the same in all three cases, and for any given multi-index $J$
in that range we may extract a common factor
\[
s \frac{(-1)^{|J|}(m - s - 1)! |J|!}{r^{|J|+1} (m - s + |J|)! J!}
d_{J + 1_{i_1}} S^{J + 1_j}
\]
from the terms in the three sums to give
\[
(m - s) - (m - s + |J| + 1) + (|J| + 1) = 0 \, .
\]

Next, we take the remaining expression and combine sums~(2) and~(3), and
also combine sums~(4) and~(7), to give
\begin{eqnarray*}
\makebox[5ex][l]{$s \, d_{i_1} P^j_{(s)} + \d^j_{i_1} P^k_{(s+1)} d_k
- s \, P^j_{(s+1)} d_{i_1}$} \\[1ex]
& = & \d^j_{i_1} \sum_{|J|=0}^{rk-1} \frac{(-1)^{|J|}
(m - s)! |J|!}{r^{|J|+1} (m - s + |J| + 1)! J!} d_{J + 1_q} S^{J + 1_q} \\
& & \phantom{skip} + \d^j_{i_1} \sum_{|J|=0}^{rk} \frac{(-1)^{|J|}
(m - s)! |J|!}{r^{|J|} (m - s + |J|)! J!} d_J S^J \\
& = & \d^j_{i_1} \sum_{|J|=0}^{rk-1} \frac{(-1)^{|J|}
(m - s)! |J|!}{r^{|J|+1} (m - s + |J| + 1)! J!} d_{J + 1_q} S^{J + 1_q} \\
& & \phantom{skip} + \d^j_{i_1} \sum_{|J|=1}^{rk} \frac{(-1)^{|J|}
(m - s)! |J|!}{r^{|J|} (m - s + |J|)! J!} d_J S^J \\
& & \phantom{skip} + \d^j_{i_1} \, ,
\end{eqnarray*}
where we have separated out the term $|J|=0$ from the second sum to give the
required $\d^j_{i_1}$. But we can rewrite this modified second sum as
\newpage
\begin{eqnarray*}
\makebox[5ex][l]{$ \displaystyle \d^j_{i_1} \sum_{|J|=1}^{rk} \frac{(-1)^{|J|}
(m - s)! |J|!}{r^{|J|} (m - s + |J|)! J!} d_J S^J$} \\[1ex]
& = & \d^j_{i_1} \sum_{|J|=0}^{rk-1} \frac{(-1)^{|J|+1}
(m - s)! (|J|+1)!}{r^{|J|+1} (m - s + |J| + 1)! (J + 1_q)!}
d_{J + 1_q} S^{J + 1_q} \times \frac{J(q) + 1}{|J|+1} \\
& = & \d^j_{i_1} \sum_{|J|=0}^{rk-1} \frac{(-1)^{|J|+1}
(m - s)! |J|!}{r^{|J|+1} (m - s + |J| + 1)! J!} d_{J + 1_q} S^{J + 1_q}
\end{eqnarray*}
where the weight $(J + 1_q)! / (|J| + 1)!$ of the multi-index $J + 1_q$ has
been replaced by the weight $J! / |J|!$ of the multi-index $J$ to take account
of the change of notation; thus the modified second sum cancels with the first
sum, and the proof is complete.
\blob

\section{The homogeneous variational bicomplexes}

In order to obtain make use of the homogeneous variational complexes, we embed
them in bicomplexes where the commuting map is the usual exterior derivative
$d : \Ors_k \to \Oras_k$. In the diagram below, we give an example of such a
bicomplex where the scalar forms are of order $k \geq 0$; it should be
understood that one or more upper rows of the diagram will have to be omitted
if $-m < k < 0$. We also write
\[
\Xi^r_{k+m} = \Orm_{k+m} / \dT \Orma_{k+m-1} \, , \qquad
\Xibar^0_{k+m} = \Oombar_{k+m} / \dT \Oomabar_{k+m-1}
\]
for the final quotient spaces in each column.

\setlength{\unitlength}{0.9pt}
\begin{center}
\begin{picture}(450,490)(0,-170)
\multiput(160,300)(80,0){3}{\vector(0,-1){40}}
\multiput(160,220)(80,0){3}{\vector(0,-1){40}}
\multiput(160,120)(80,0){3}{\makebox(0,0){$\vdots$}}
\multiput(160,60)(80,0){3}{\vector(0,-1){40}}
\multiput(160,-20)(80,0){3}{\vector(0,-1){40}}
\multiput(160,-100)(80,0){3}{\vector(0,-1){40}}
\multiput(105,240)(0,-80){2}{\vector(1,0){30}}
\put(110,80){\vector(1,0){25}}
\multiput(105,0)(0,-80){2}{\vector(1,0){30}}
\multiput(180,240)(0,-80){5}{\vector(1,0){35}}
\put(180,-20){\vector(1,-1){40}}
\multiput(260,240)(0,-80){5}{\vector(1,0){35}}
\multiput(360,240)(0,-80){5}{\makebox(0,0){$\cdots$}}
\multiput(80,240)(0,-80){5}{\makebox(0,0){$0$}}
\put(160,320){\makebox(0,0){$0$}}
\multiput(240,320)(80,0){2}{\makebox(0,0){$0$}}
\put(160,240){\makebox(0,0){$\Ooobar_k$}}
\put(240,240){\makebox(0,0){$\Oao_k$}}
\put(320,240){\makebox(0,0){$\Obo_k$}}
\put(160,160){\makebox(0,0){$\Ooabar_{k+1}$}}
\put(240,160){\makebox(0,0){$\Oaa_{k+1}$}}
\put(320,160){\makebox(0,0){$\Oba_{k+1}$}}
\put(160,80){\makebox(0,0){$\Oomabar_{k+m-1}$}}
\put(240,80){\makebox(0,0){$\Oama_{k+m-1}$}}
\put(320,80){\makebox(0,0){$\Obma_{k+m-1}$}}
\put(160,0){\makebox(0,0){$\Oombar_{k+m}$}}
\put(240,0){\makebox(0,0){$\Oam_{k+m}$}}
\put(320,0){\makebox(0,0){$\Obm_{k+m}$}}
\put(160,-80){\makebox(0,0){$\Xibar^0_{k+m}$}}
\put(240,-80){\makebox(0,0){$\Xi^1_{k+m}$}}
\put(320,-80){\makebox(0,0){$\Xi^2_{k+m}$}}
\multiput(160,-160)(80,0){3}{\makebox(0,0){$0$}}
\multiput(200,10)(0,80){4}{\makebox(0,0)[b]{$\scriptstyle d$}}
\multiput(280,10)(0,80){4}{\makebox(0,0)[b]{$\scriptstyle d$}}
\multiput(200,-70)(80,0){2}{\makebox(0,0)[b]{$\scriptstyle \d$}}
\multiput(155,200)(80,0){3}{\makebox(0,0)[r]{$\scriptstyle\dT$}}
\multiput(155,40)(80,0){3}{\makebox(0,0)[r]{$\scriptstyle\dT$}}
\put(155,-40){\makebox(0,0)[r]{$\scriptstyle p_0$}}
\put(235,-40){\makebox(0,0)[r]{$\scriptstyle p_1$}}
\put(315,-40){\makebox(0,0)[r]{$\scriptstyle p_2$}}
\end{picture}
\end{center}
\setlength{\unitlength}{1pt}

\setlength{\unitlength}{0.9pt}
\begin{center}
\begin{picture}(400,600)(0,-180)
\multiput(-40,240)(0,-80){5}{\makebox(0,0){$\cdots$}}
\put(0,300){\vector(0,-1){40}}
\put(0,220){\vector(0,-1){40}}
\multiput(0,120)(80,0){2}{\makebox(0,0){$\vdots$}}
\multiput(0,60)(80,0){3}{\vector(0,-1){40}}
\multiput(0,-20)(80,0){3}{\vector(0,-1){40}}
\multiput(0,-100)(80,0){3}{\vector(0,-1){40}}
\multiput(20,240)(0,-80){2}{\vector(1,0){35}}
\put(23,80){\vector(1,0){25}}
\multiput(20,0)(0,-80){2}{\vector(1,0){35}}
\multiput(120,160)(0,-80){4}{\makebox(0,0){$\cdots$}}
\put(190,80){\vector(1,0){35}}
\put(185,0){\vector(1,0){22}}
\put(180,-60){\vector(1,1){40}}
\put(280,0){\makebox(0,0){$\cdots$}}
\put(340,0){\vector(1,0){40}}
\put(0,240){\makebox(0,0){$\Omega^{N_0,0}_k$}}
\put(0,160){\makebox(0,0){$\Omega^{N_0,1}_{k+1}$}}
\put(0,80){\makebox(0,0){$\Omega^{N_0,m-1}_{k+m-1}$}}
\put(0,0){\makebox(0,0){$\Omega^{N_0,m}_{k+m}$}}
\put(80,240){\makebox(0,0){$0$}}
\put(80,160){\makebox(0,0){$\Omega^{N_0+1,1}_{k+1}$}}
\put(80,80){\makebox(0,0){$\Omega^{N_0+1,m-1}_{k+m-1}$}}
\put(80,0){\makebox(0,0){$\Omega^{N_0+1,m}_{k+m}$}}
\put(0,320){\makebox(0,0){$0$}}
\put(160,80){\makebox(0,0){$\Omega^{N_{m-1},m-1}_{k+m-1}$}}
\put(240,80){\makebox(0,0){$0$}}
\put(160,0){\makebox(0,0){$\Omega^{N_{m-1},m}_{k+m}$}}
\put(240,0){\makebox(0,0){$\Omega^{N_{m-1}+1,m}_{k+m}$}}
\put(320,0){\makebox(0,0){$\Omega^{N_m,m}_{k+m}$}}
\put(0,-80){\makebox(0,0){$\Xi^{N_0}_{k+m}$}}
\put(80,-80){\makebox(0,0){$\Xi^{N_0+1}_{k+m}$}}
\put(160,-80){\makebox(0,0){$\Xi^{N_{m-1}}_{k+m}$}}
\put(400,0){\makebox(0,0){$0$}}
\multiput(0,-160)(80,0){3}{\makebox(0,0){$0$}}
\multiput(40,10)(0,80){4}{\makebox(0,0)[b]{$\scriptstyle d$}}
\multiput(200,10)(0,80){2}{\makebox(0,0)[b]{$\scriptstyle d$}}
\put(360,10){\makebox(0,0)[b]{$\scriptstyle d$}}
\put(40,-70){\makebox(0,0)[b]{$\scriptstyle \d$}}
\put(-5,200){\makebox(0,0)[r]{$\scriptstyle\dT$}}
\multiput(-5,40)(80,0){3}{\makebox(0,0)[r]{$\scriptstyle\dT$}}
\put(-5,-40){\makebox(0,0)[r]{$\scriptstyle p_{N_0}$}}
\put(75,-40){\makebox(0,0)[r]{$\scriptstyle p_{N_0+1}$}}
\put(155,-40){\makebox(0,0)[r]{$\scriptstyle p_{N_{m-1}}$}}
\end{picture}
\end{center}
\setlength{\unitlength}{1pt}
On the right-hand side of the diagram, the numbers $N_s$ represent the
dimensions of the manifolds $\Fmks E$. Note that $d$ and $\dT$ commute
rather than anti-commute.

We have proved that the columns of this diagram (apart from the first) are
exact modulo pullbacks, and we know that the rows (apart from the last) are
locally exact by the standard property of the exterior derivative:\ recall that
we have replaced the usual zeroth term of the de Rham sequence with its
quotient by the constants. Local exactness elsewhere (again, possibly modulo
pull-backs) comes from diagram chasing. The following arguments are based on
those in~\cite{Tul}, although given there in a slightly different context. The
assertions hold only locally, so a symbol such as $\Ors$ should be
interpreted here as a sheaf of germs of vector-valued forms, rather than as a
space of globally-defined forms. (We shall not be specific about the orders
of the forms here, and we shall again omit the pull-back maps, as the
formul\ae\ for the orders become increasingly complex and obscure the main
thrust of the argument.) Some modifications will be needed in the case of
zeroth-order forms taking their values in $\bw^s \R^{m*}$ for $s \geq 1$, but
the nature of these modifications should be clear. The gradual shortening of
the columns on the right-hand side of the bicomplex has no effect on our
argument, as we always increase the order when moving upwards.

\textbf{Lemma 5}
{\em If $1 \leq r \leq N_0 - 1$ and $1 \leq s \leq m - 1$ then
\[
\ker \left( \dT d: \Oro \to \Oraa \right) \subset
\im \left( d : \OrAo \to \Oro \right)
\]
and
\[
\ker \left( \dT d : \Ors \to \Orasa \right) \subset
\im \left( d : \OrAs \to \Ors \right)
+ \im \left( \dT : \OrsA \to \Ors \right) \, .
\]
}

\textbf{Proof}
Suppose first that
\[
\Phi \in \ker \left( \dT d: \Oro \to \Oraa \right) \, ,
\]
so that $\dT d\phi = 0$. Then $d\phi = 0$ because $\dT : \Oro \to \Ora$
is injective; thus
\[
\Phi \in \im \left( d : \OrAo \to \Oro \right) \, .
\]

The second assertion follows from the first by induction on $s$. Suppose
that
\[
\Phi \in \ker \left( \dT d: \Ors \to \Orasa \right) \, ,
\]
so that
\[
d\Phi \in \ker \left( \dT: \Oras \to \Orasa \right) \, .
\]
Put $\Phi_0 = Pd\Phi \in \OrasA$, so that $\dT \Phi_0 = d\Phi$. Then
\[
\dT d\Phi_0 = d\dT \Phi_0 = d^2 \Phi = 0 \, ,
\]
so that
\[
\Phi_0 \in \ker \left( \dT d: \OrasA \to \Orbs \right) \, .
\]
Thus by the induction hypothesis we may write
\[
\Phi_0 = d\Phi_1 + \dT \Phi_2
\]
where
\[
\Phi_1 \in \OrsA \, , \qquad \Phi_2 \in \OrasB  \, .
\]
Now consider $\Phi - \dT \Phi_1$. We have
\begin{eqnarray*}
d \left( \Phi - \dT \Phi_1 \right)
& = & \dT \Phi_0  - \dT \left( \Phi_0 - \dT \Phi_2 \right) \\
& = & 0 \, ,
\end{eqnarray*}
so that
\[
\Phi - \dT \Phi_1 = d\Psi
\]
for some $\Psi \in \OrAs$ and the result follows.
\blob

\textbf{Lemma 6}
{\em If $1 \leq s \leq m - 1$ then
\[
\ker \left( \dT d: \Ooobar \to \Oaa \right) = 0
\]
and
\[
\ker \left( \dT d : \Oos \to \Oasa \right)
\subset \im \left( \dT : \OosAbar \to \Oosbar \right) \, .
\]
}

\textbf{Proof}
Take the results of Lemma~5 with $r = 1$ and apply the same method of proof
once again.
\blob

\textbf{Corollary}
{\em
The left-hand column of each homogeneous variational bicomplex is locally
exact.
}
\blob

It should now be clear why we have chosen to use the spaces $\Oosbar_{k+s}$
rather than $\Oos_{k+s}$:\ the constant vector-valued functions are
certainly $\dT$-closed, but are not locally $\dT$-exact, even to within
pullback. This is rather different from the behaviour of the left-hand column
of the variational bicomplex on jet bundles, which {\em is}\/ locally
$\dh$-exact on forms rather than classes of forms.

\textbf{Lemma 7}
{\em
If $1 \leq r \leq N_{m-1}$ then
\[
\ker \left( \d p_r : \Orm \to \Xi^{r+1} \right)
\subset \im \left( d : \OrAm \to \Orm \right)
+ \im \left( \dT : \Orma \to \Orm \right) \, ,
\]
and in addition
\[
\ker \left( \d p_0 : \Oombar \to \Xi^1 \right)
\subset \im \left( \dT : \Oomabar \to \Oombar \right) \, .
\]
}

\textbf{Proof}
Take the results of Lemma~5 with $s = m$ and apply the same method of proof
once again.
\blob

\textbf{Corollary}
{\em
The bottom row of each homogeneous variational bicomplex is locally exact.
}
\blob

\section{Homogeneous Lagrangians and the edge sequence}

The {\em edge sequence}\/ of a homogeneous variational bicomplex, as given
above, is
\[
0 \to \Ooobar_k \to \Ooabar_{k+1} \to \ldots \to \Oosbar_{k+s} \to
\Oosabar_{k+s+1} \to \ldots \to \hspace{4em}
\]
\\[-6ex]
\[
\hspace{4em} \to \Oombar_{k+m} \to
\Xi^1_{k+m} \to \ldots \to \Xi^{N_{m-1}}_{k+m} \to \Omega^{N_{m-1}+1,m}_{k+m}
\to \ldots \to \Omega^{N_m,m}_{k+m} \to 0 \, ,
\]
and we have seen that this is locally exact, modulo pull-backs. This sequence
is particularly significant for homogeneous variational problems. In this final
section, it is convenient to relabel the order of the various spaces so that
we have $\Oombar_k$ rather than $\Oombar_{k+m}$; if $k < 0$ we shall, of course,
need to omit the first $-k$ non-zero terms of the edge sequence.

A homogeneous $k$-th order Lagrangian may be considered as a function $L$ on an
$m$-frame bundle $\Fmk E$ satisfying the conditions
\[
\D^i_j(L) = \d^i_j L \, ,
\qquad \D^I_j(L) = 0 \quad \mbox{for} \quad |I| > 1 \, .
\]
where $\D^I_j = S^I(T_j)$ (note that, although $T_j$ is a vector field along,
in this case, the projection $\taumkkaE$, the contraction is unambiguous). The
Hilbert forms $\vt^i$ of the Lagrangian are defined by
\[
\vt^i = P^i_{(1)} dL \, ,
\]
and the Euler-Lagrange form $\ve$ by
\[
\ve = dL - d_i \vt^i
\]
so that, in coordinates,
\[
\ve = \sum_{|I|=0}^k (-1)^{|I|} d_I \left( \pd{L}{u^\a_I} \right)
\]
(see, for example,~\cite{CS2} for further details of these constructions).

Now the function $L$ may be considered as the vector valued function
$\L = L \, d^m t$, where $d^m t = dt^1 \wedge \ldots \wedge dt^m$ is the
canonical volume element on $\R^m$, and hence we may determine whether an
element $\L \in \Oom_k E$ is homogeneous. An element of the quotient space
$\Oombar_k$ will be said to be homogeneous if it has a homogeneous
representative; such a representative must be unique as the only homogeneous
constant function is zero. The individual Hilbert forms may be combined into
a single vector-valued form
\[
\Theta_1 = \vt^i \otimes \left( \vf{t_i} \hook d^m t \right)
= PdL \, ,
\]
and then the formula for the vector-valued version of the Euler-Lagrange form
becomes
\[
\ve \otimes d^m t = d\L - \dT \Theta_1 \in \Oam_{2k} \, .
\]
Indeed, the pseudo-homotopy formula for $\dT$ tells us that the map
$\Xi^r_k \to \Orm_k$ given by
\[
p_r(\Phi) \mapsto \Phi - \dT P\Phi
\]
is globally well-defined for $r \geq 1$, and gives a canonical global
representative for each class in $\Xi^r_k$; the Euler-Lagrange form
$\ve \otimes d^m t$ is the canonical representative of the class $p_1(dL)$.
In a similar way, if $\Phi \in \Oam_{2k}$ then the map
$\H : \Oam_{2k} \to \Obm_{4k}$ given by
\[
\H(\Phi) = d\Phi - \dT P\Phi
\]
is called the {\em Helmholtz-Sonin map}\/, and if $\H(\Phi) = 0$ then
local exactness of the bicomplex shows that $\Phi$ must locally be an
Euler-Lagrange form.

There are also questions that may be answered by studying parts of the
bicomplex away from the edge sequence. One such question involves the existence
of a scalar $m$-form $\Theta_m$ corresponding to a Lagrangian $\L$, having the
property that $\Theta_m$ is closed precisely when $\L$ is a null Lagrangian
(that is, when $\ve \otimes d^m t = 0$). For a first-order Lagrangian, such an
$m$-form is given by
\[
\Theta_m = S^1 d S^2 d \ldots S^m dL
\]
(see~\cite{CS3}); using the language of vector-valued forms described above, we
may write this as
\[
\Theta_m = (Pd)^m \L \in \Omo_1
\]
where the formula $(Pd)^m \L$ gives in principle a form of higher order, but
the result is projectable to the first-order frame bundle. One might ask
whether a similar result is true for $k$-th order Lagrangians with $k > 1$.
Some preliminary work~\cite{SauFund} shows that this is indeed the case for
second-order Lagrangians in two independent variables, where $\Theta_2$ is
projectable to the fourth-order frame bundle, and investigation of the more
general problem continues.

In summary, therefore, we conclude that the homogeneous variational bicomplexes
may be used to analyse homogeneous variational problems on $m$-frame bundles in
the same way that the variational bicomplex and finite-order variational
sequence are used for variational problems on jet bundles.

\end{document}